\documentclass[11pt,fleqn]{amsart2000}
\newtheorem{theorem}{Theorem}[section]

\newtheorem{corollary}[theorem]{Corollary}

\begin{document}
\setlength{\unitlength}{0.01in}
\linethickness{0.01in}
\begin{center}
\begin{picture}(474,66)(0,0)
\multiput(0,66)(1,0){40}{\line(0,-1){24}}
\multiput(43,65)(1,-1){24}{\line(0,-1){40}}
\multiput(1,39)(1,-1){40}{\line(1,0){24}}
\multiput(70,2)(1,1){24}{\line(0,1){40}}
\multiput(72,0)(1,1){24}{\line(1,0){40}}
\multiput(97,66)(1,0){40}{\line(0,-1){40}}
\put(143,66){\makebox(0,0)[tl]{\footnotesize Proceedings of the Ninth Prague Topological Symposium}}
\put(143,50){\makebox(0,0)[tl]{\footnotesize Contributed papers from the symposium held in}}
\put(143,34){\makebox(0,0)[tl]{\footnotesize Prague, Czech Republic, August 19--25, 2001}}
\end{picture}
\end{center}
\vspace{0.25in}
\setcounter{page}{147}
\title[Fell-continuous selections]{Fell-continuous selections and
topologically well-orderable spaces II}
\author{Valentin Gutev}
\address{School of Mathematical and Statistical Sciences,
Faculty of Science, University of Natal, King George V Avenue,
Durban 4041, South Africa}
\email{gutev@nu.ac.za}
\thanks{Valentin Gutev,
{\em Fell-continuous selections and topologically well-orderable spaces II},
Proceedings of the Ninth Prague Topological Symposium, (Prague, 2001),
pp.~147--153, Topology Atlas, Toronto, 2002}
\begin{abstract}
The present paper improves a result of \cite{gutev:00d} by
showing that a space $X$ is topologically well-orderable if and only
if there exists a selection for $\mathcal{F}_2(X)$ which is
continuous with respect to the Fell topology on $\mathcal{F}_2(X)$.
In particular, this implies that $\mathcal{F}(X)$ has a
Fell-continuous selection if and only if $\mathcal{F}_2(X)$ has a
Fell-continuous selection.
\end{abstract}
\maketitle

\section{Introduction}

Let $X$ be a topological space, and let $\mathcal{F}(X)$ be the family
of all non-empty closed subsets of $X$. Also, let $\tau$ be a topology
on $\mathcal{F}(X)$ and $\mathcal{D}\subset \mathcal{F}(X)$. A map
$f:\mathcal{D}\to X$ is a \emph{selection} for $\mathcal{D}$ if
$f(S)\in S$ for every $S\in\mathcal{D}$. A map $f:\mathcal{D}\to
X$ is a \emph{$\tau$-continuous} selection for $\mathcal{D}$ if it
is a selection for $\mathcal{D}$ which is continuous with respect
to the relative topology $\tau$ on $\mathcal{D}$ as a subspace of
$\mathcal{F}(X)$.

Two topologies on $\mathcal{F}(X)$ will play the most important role
in this paper. The first one is the \emph{Vietoris topology} $\tau_V$
which is generated by all collections of the form
\[
\langle\mathcal{V}\rangle = 
\left\{S\in \mathcal{F}(X) : 
S\cap V\ne \emptyset,\ 
V\in \mathcal{V},\ 
\mbox{and}\ 
S\subset \bigcup \mathcal{V}\right\},
\]
where $\mathcal{V}$ runs over the finite families of open subsets of
$X$. The other one is the \emph{Fell topology} $\tau_F$ which is defined
by all basic Vietoris neighbourhood $\langle \mathcal{V}\rangle$ with
the property that $X\setminus \bigcup\mathcal{V}$ is compact.

Finally, let us recall that a space $X$ is \emph{topologically
well-orderable} (see Engelking, Heath and Michael
\cite{engelking-heath-michael:68}) if there exists a linear order
``$\prec$'' on $X$ such that $X$ is a linear ordered topological space
with respect to $\prec$, and every non-empty closed subset of $X$
has a $\prec$-minimal element.

Recently, the topologically well-orderable spaces were characterized
in \cite[Theorem 1.3]{gutev:00d} by means of Fell-continuous
selections for their hyperspaces of non-empty closed subsets.

\begin{theorem}[\cite{gutev:00d}]\label{th:Fell-1}
A Hausdorff space $X$ is topologically well-orderable if and only if
$\mathcal{F}(X)$ has a $\tau_F$-continuous selection.
\end{theorem}

In the present paper, we improve Theorem \ref{th:Fell-1} by showing
that one may use $\tau_F$-continuous selections only for the subset
$\mathcal{F}_2(X)=\{S\in \mathcal{F}(X): |S|\le 2\}$ of
$\mathcal{F}(X)$. Namely, the following theorem will be proven. 

\begin{theorem}\label{th:Fell-2}
A Hausdorff space $X$ is topologically well-orderable if and only if
$\mathcal{F}_2(X)$ has a $\tau_F$-continuous selection.
\end{theorem}

About related results for Vietoris-continuous selections, the interested
reader is referred to van Mill and Wattel \cite{mill-wattel:81}.

Theorem \ref{th:Fell-2} is interesting also from another point of
view. According to Theorem \ref{th:Fell-1}, it implies the following
result which may have an independent interest.

\begin{corollary}\label{cr:Fell-3}
If $X$ is a Hausdorff space, then $\mathcal{F}(X)$ has a 
$\tau_F$-contin\-uous selection if and only if $\mathcal{F}_2(X)$ has a
$\tau_F$-continuous selection.
\end{corollary}

A word should be said also about the proof of Theorem \ref{th:Fell-2}.
In general, it is based on the proof of Theorem \ref{th:Fell-1} stated
in \cite{gutev:00d}, and is separated in a few different steps which
are natural generalizations of the corresponding ones given in
\cite{gutev:00d}. In fact, the paper demonstrates that all statements
of \cite{gutev:00d} remain true if $\mathcal{F}(X)$ is replaced by
$\mathcal{F}_{2}(X)$. Related to this, the interested reader may
consult an alternative proof of Theorem \ref{th:Fell-2} given in
\cite{artico-marconi:00} and based again on the scheme in
\cite{gutev:00d}.

\section{A reduction to locally compact spaces}

In the sequel, all spaces are assumed to be at least Hausdorff.

In this section, we prove the following generalization of
\cite[Theorem 2.1]{gutev:00d}.

\begin{theorem}\label{th:step1}
Let $X$ be a space such that $\mathcal{F}_{2}(X)$ has a 
$\tau_F$-continuous selection. 
Then $X$ is locally compact.
\end{theorem}

\begin{proof}
We follow the proof of \cite[Theorem 2.1]{gutev:00d}. 
Namely, let $f$ be a $\tau_F$-continuous selection for $\mathcal{F}_2(X)$
and suppose, if possible, that $X$ is not locally compact. 
Hence, there exists a point $p\in X$ such that $\overline{V}$ is not 
compact for every neighbourhood $V$ of $p$ in $X$. 
Claim that there exists a point $q\in X$ such that 
\begin{equation}\label{eq:point-q}
q\ne p\ \mbox{and}\ f(\{p,q\})=p.
\end{equation}
To this purpose, note that there exists $F\in\mathcal{F}(X)$ such that $F$
is not compact and $p\notin F$. 
Then, $f^{-1}(X\setminus F)$ is a $\tau_F$-neighbourhood of $\{p\}$ in
$\mathcal{F}_2(X)$, so there exists a finite family $\mathcal{W}$ of open
subsets of $X$ such that $X\setminus \bigcup \mathcal{W}$ is compact and
\[
\{p\}\in\langle \mathcal{W}\rangle\cap\mathcal{F}_2(X)
\subset f^{-1}(X\setminus F).
\]
Then, $F\cap W\neq\emptyset$ for some $W\in \mathcal{W}$ because $F$ is
not compact. 
Therefore, there exists a point 
$q\in F\cap \left(\bigcup \mathcal{W}\right)$. 
This $q$ is as required.

Let $q$ be as in (\ref{eq:point-q}). 
Since $X$ is Hausdorff, $f(\{q\})\neq f(\{p,q\})$, and $f$ is
$\tau_F$-continuous, there now exist two finite families $\mathcal{U}$
and $\mathcal{V}$ of open subsets of $X$ such that 
$X\setminus \bigcup \mathcal{U}$ is compact, 
$\{q\}\in \langle \mathcal{U}\rangle$, 
$\{p,q\}\in \langle \mathcal{V}\rangle$, and 
$\langle \mathcal{U}\rangle\cap \langle \mathcal{V}\rangle=\emptyset$. 
Then,
\begin{equation}\label{eq:point}
p\in V_p = \bigcap\{V\in \mathcal{V}: p\in V\} \subset 
X\setminus \bigcup \mathcal{U}.
\end{equation}
Indeed, suppose there is a point $\ell\in V_p\cap \left(\bigcup
\mathcal{U} \right)$. Then, $\{\ell,q\}\in\langle
\mathcal{U}\rangle$ because $\{q\} \in \langle \mathcal{U}\rangle$.
However, we also get that $\{\ell,q\}\in\langle \mathcal{V}\rangle$
because $q\notin V$ for some $V\in \mathcal{V}$ implies $p\in V$,
hence $\ell\in V_p\subset V$. Thus, we finally get that
$\{\ell,q\}\in\langle \mathcal{U}\rangle\cap\langle
\mathcal{V}\rangle$ which is impossible. So, (\ref{eq:point})
holds as well.

To finish the proof, it remains to observe that this contradicts the
choice of $p$. Namely $V_p$ becomes a neighbourhood of $p$ which, by
(\ref{eq:point}), has a compact closure because $X\setminus \bigcup
\mathcal{U}$ is compact.
\end{proof}

\section{A reduction to compact spaces}

For a locally compact space $X$ we will use $\alpha X$ to denote the one
point compactification of $X$. For a non-compact locally compact $X$ let
us agree to denote by $\alpha$ the point of the singleton 
$\alpha X\setminus X$.

In what follows, to every family $\mathcal{D}\subset \mathcal{F}(X)$ we
associate a family $\alpha(\mathcal{D})\subset \mathcal{F}(\alpha X)$
defined by
\[
\alpha(\mathcal{D}) = 
\{S\in \mathcal{F}(\alpha X): S\cap X\in \mathcal{D}\cup \{\emptyset\}\}.
\]

The following extension theorem was actually proven in \cite[Theorem 
3.1]{gutev:00d}.

\begin{theorem}\label{th:step2}
Let $X$ be a locally compact non-compact space $X$, and
$\mathcal{D}\subset \mathcal{F}(X)$. Then, $\mathcal{D}$ has a
$\tau_F$-continuous selection if and only if $\alpha(\mathcal{D})$
has a $\tau_V$-continuous selection $g$ such that
$g^{-1}(\alpha)=\{\{\alpha\}\}$.
\end{theorem}

\begin{proof} Just the same proof as in \cite[Theorem
3.1]{gutev:00d} works. Namely, if $f$ is a $\tau_{F}$-continuous
selection for $\mathcal{D}$, we may define a selection $g$ for
$\alpha(\mathcal{D})$ by $g(S)=f(S\cap X)$ if $S\cap X\neq\emptyset$
and $g(S)=\alpha$ otherwise, where $S\in \alpha(\mathcal{D})$.
Clearly $g^{-1}(\alpha)=\{\{\alpha\}\}$ and, as shown in
\cite[Theorem 3.1]{gutev:00d}, $g$ is $\tau_V$-continuous. If now
$g$ is a $\tau_{V}$-continuous selection for $\alpha(\mathcal{D})$,
with $g^{-1}(\alpha)=\{\{\alpha\}\}$, then $g(S\cup\{\alpha\})\in S$
for every $S\in \mathcal{D}$, so we may define a selection $f$ for
$\mathcal{D}$ by $f(S)=g\left(S\cup\{\alpha\}\right)$,
$S\in\mathcal{D}$. The verification that $f$ is $\tau_F$-continuous
was done in \cite[Theorem 3.1]{gutev:00d}.
\end{proof}

\section{Special selections and connected sets}

In what follows, to every selection $f:\mathcal{F}_2(X)\to X$ we
associate an order-like relation ``$\prec_f$'' on $X$ (see Michael
\cite{michael:51}) defined for $x\ne y$ by
\[
x_1\prec_f x_2\
\mbox{iff}\
f(\{x_1,x_2\})=x_1.
\]
Further, we will need also the following $\prec_{f}$-intervals:
\[
(x,+\infty)_{\prec_f}=\{z\in X:x\prec_f z\}
\]
and
\[
[x,+\infty)_{\prec_f}=\{z\in X:x\preceq_f z\}.
\]
Now, we provide the generalization of \cite[Theorem 4.1]{gutev:00d}
for the case of $\mathcal{F}_{2}(X)$.

\begin{theorem}\label{th:connected}
Let $X$ be a space, $a\in X$, and let $A\in \mathcal{F}(X)$ be a connected
set such that $|A|>1$ and $a\in A\cap \overline{X\setminus A}$. 
Also, let $f:\mathcal{F}_{2}(X)\to X$ be a $\tau_V$-continuous selection
for $\mathcal{F}_2(X)$. 
Then, $f^{-1}(a)\neq\{\{a\}\}$.
\end{theorem}

\begin{proof}
Suppose, on the contrary, that $f^{-1}(a)=\{\{a\}\}$. By
hypothesis, there exists a point $b\in A$, with $b\ne a$. Since $f$
is $\tau_{V}$-continuous, $f(\{a,b\})=b$ and $a\in
\overline{X\setminus A}$, we can find a point $c\in X\setminus A$
such that $f(\{b,c\})=b$. Then, $B=A\cap (c,+\infty)_{\prec_{f}}$ is
a clopen subset of $A$ because $B=A\cap [c,+\infty)_{\prec_{f}}$,
see \cite{michael:51}. However, this is impossible because $b\in
A\setminus B$, while $a\in B$.
\end{proof}

\section{A further result about special selections}

Following \cite{gutev:00d}, we shall say that a point $a\in X$ is a
\emph{partition} of $X$ if there are open subset $L,R\subset
X\setminus \{a\}$ such that $\overline{L}\cap\overline{R}=\{a\}$ and
$L\cap R=\emptyset$.

We finalize the preparation for the proof of Theorem \ref{th:Fell-2} with
the following result about special Vietoris continuous selections and
partitions which generalizes \cite[Theorem 5.1]{gutev:00d}.

\begin{theorem}\label{th:partit}
Let $X$ be a compact space, $f$ a $\tau_V$-continuous selection
for $\mathcal{F}_2(X)$, and let $a\in X$ be a partition of $X$ such
that $f^{-1}(a)=\{\{a\}\}$. Then, $X$ is first countable at $a$.
\end{theorem}

\begin{proof}
By definition, there are open sets $L,R\subset X\setminus\{a\}$ such that
$\overline{L}\cap\overline{R}=\{a\}$ and $L\cap R=\emptyset$.
Hence, both $L$ and $R$ are non-empty. 
Take a point $\ell_0\in L$.
Then, by hypothesis, $f(\{\ell_0,a\})=\ell_0$. 
Since $f$ is $\tau_V$-continuous, this implies the existence of a
neighbourhood $L_0\subset L$ of $\ell_0$ and a neighbourhood $V_0$ of $a$
such that
\[
L_0\cap V_0 = \emptyset
\mbox{ and }
f(\langle \{L_0,V_0\}\rangle\cap\mathcal{F}_2(X))\subset L_0.
\]
Since $a\in\overline{R}$, there exists a point $r_0\in V_0\cap R$. 
Observe that $f(\{a,r_0\})=r_0\in V_0$. 
Hence, just like before, we may find a neighbourhood 
$R_0\subset R\cap V_0$ of $r_0$ and a neighbourhood $W_0\subset V_0$ of
$a$ such that
\[
R_0\cap W_0 = \emptyset
\mbox{ and }
f(\langle \{R_0,W_0\}\rangle\cap\mathcal{F}_2(X))\subset R_0.
\]
Thus, by induction, we may construct a sequence $\{\ell_n:n<\omega\}$ of
points of $L$, a sequence $\{r_n:n<\omega\}$ of points of $R$, and open
sets $L_n,V_n,R_n,W_n\subset X$ such that
\begin{equation}\label{eq:ln}
\begin{array}{l}
\ell_n\in L_n,\\
a\in V_n,\\
L_n\cap V_n=\emptyset \mbox{ and}\\
f(\langle \{L_n,V_n\}\rangle \cap\mathcal{F}_2(X))\subset L_n,
\end{array}
\end{equation}
\begin{equation}\label{eq:rn}
\begin{array}{l}
r_n\in R_n,\\
a\in W_n,\\
R_n\cap W_n=\emptyset \mbox{ and}\\
f(\langle \{R_n,W_n\} \rangle \cap\mathcal{F}_2(X))\subset R_n,
\end{array}
\end{equation}
and
\begin{equation}\label{eq:lr}
\begin{array}{l}
{V_{n+1}}\subset W_n\subset V_n,\\
L_{n+1}\subset L\cap W_n \mbox{ and}\\
R_{n}\subset R\cap V_n.
\end{array}
\end{equation}
Since $X$ is compact, $\{\ell_{n}:n<\omega\}$ has a cluster point
$\ell$, and $\{r_{n}:n<\omega\}$ has a cluster point $r$. We claim
that $\ell=r$. Indeed, suppose for instance that $\ell\prec_{f} r$
(the case $r\prec_{f} \ell$ is symmetric). Then, there are disjoint
open sets $U_{\ell}$ and $U_{r}$ such that $\ell\in U_{\ell}$, $r\in
U_{r}$, and $x\prec_{f} y$ for every $x\in U_{\ell}$ and $y\in
U_{r}$, see \cite{gutev-nogura:01a}. Next, take $\ell_{n}\in
U_{\ell}$ and $r_{m}\in U_{r}$ such that $n>m$. Then, we have
$\ell_{n}\prec_{f} r_{m}$. However, by (\ref{eq:ln}), (\ref{eq:rn})
and (\ref{eq:lr}), we get that $\{r_{m},\ell_{n}\}\in
\langle\{R_{m},W_{m}\}\rangle\cap\mathcal{F}_2(X)$, and therefore
$f(\{r_{m},\ell_{n}\})=r_{m}$. This is clearly impossible, so
$\ell=r$.

Having already established this, let us observe that $b=\ell=r$
implies $b\in \overline{L}\cap \overline{R}$ because $\ell\in
\overline{L}$ and $r\in\overline{R}$. However, $\overline{L}\cap
\overline{R}=\{a\}$ which finally implies that $b=a$.

We are now ready to prove that, for instance, $\{W_n:n<\omega\}$ is
a local base at $a$. To this end, suppose if possible that this
fails. Hence, there exists an open neighbourhood $U$ of $a$ such that
$W_n\setminus U\ne\emptyset$ for every $n<\omega$. Next, whenever
$n<\omega$, take a point $t_n\in W_n\setminus U$. Since $X$ is
compact, $\{t_n:n<\omega\}$ has a cluster point $t\not\in U$. Then,
$t\prec_{f} a$ and, as before, we may find disjoint open sets
$U_{t}$ and $U_{a}$ such that $t\in U_{t}$, $a\in U_{a}$, and
$x\prec_{f} y$ for every $x\in U_{t}$ and $y\in U_{a}$. Next, take
$t_{n}\in U_{t}$ and $r_{m}\in U_{a}$ such that $n>m$. Then,
$t_{n}\prec_{f} r_{m}$, while, by (\ref{eq:rn}) and (\ref{eq:lr}),
$r_{m}\prec_{f} t_{n}$ because $\{r_{m},t_{n}\}\in
\langle\{R_{m},W_{m}\}\rangle\cap\mathcal{F}_2(X)$.
The contradiction so obtained
completes the proof.~\end{proof}

\section{Proof of Theorem \ref{th:Fell-2}}

In case $X$ is a topologically well-orderable space, we may use
Theorem \ref{th:Fell-1}.

Suppose that $\mathcal{F}_2(X)$ has a $\tau_F$-continuous selection.
If $X$ is compact, then Theorem~\ref{th:Fell-2} is, in fact, a result
of van Mill and Wattel \cite{mill-wattel:81}. Let $X$ be non-compact.
By Theorem \ref{th:step1}, $X$ is locally compact. Then, by Theorem
\ref{th:step2}, $\mathcal{F}_2(\alpha X)$ has a $\tau_V$-continuous
selection $f$ such that $f^{-1}(\alpha)=\{\{\alpha\}\}$. Relying once
again on the result of \cite{mill-wattel:81}, $\alpha X$ is a linear
ordered topological space with respect to some linear order ``$<$'' on
$\alpha X$. It now suffices to show that there exists a compatible
(with the topology of $\alpha X$) linear order ``$\prec$'' on $\alpha
X$ such that $\alpha$ is either the first or the last element of
$\alpha X$, see \cite[Lemma 4.1]{engelking-heath-michael:68}. We show
this following precisely the proof of Theorem \ref{th:Fell-1} in
\cite{gutev:00d}. Namely, let
\[
L = \{x\in \alpha X:x<\alpha\}\
\mbox{and}\
R = \{x\in \alpha X:\alpha< x\}.
\]
Note that $L,R\subset \alpha X\setminus \{\alpha\}=X$ are open subsets
of $\alpha X$. In case one of these sets is also closed, the desired
linear order ``$\prec$'' on $\alpha X$ can be defined by exchanging
the places of $L$ and $R$. Namely, by letting for $x,y\in \alpha X$
that $x\prec y$ if and only if
\[
\begin{array}{l}
x,y\in \overline{L}\ \mbox{and}\ x<y,\ \mbox{or}\\
x,y\in \overline{R}\ \mbox{and}\ x<y,\ \mbox{or}\\
x\in \overline{R}\ \mbox{and}\ y\in \overline{L}.
\end{array}
\]
Finally, let us consider the case 
$\overline{L}\cap \overline{R}=\{\alpha\}$. 
Then, $\alpha$ is a partition of $\alpha X$.
Hence, by Theorem \ref{th:partit}, $\alpha X$ is first countable at 
$\alpha$. 
Let $\mathcal{C}[\alpha]$ be the connected component of $\alpha$ in
$\alpha X$. 
Since $f^{-1}(\alpha)=\{\{\alpha\}\}$, it now follows from Theorem
\ref{th:connected} that $\mathcal{C}[\alpha]=\{\alpha\}$.
Indeed, $\mathcal{C}'=\mathcal{C}[\alpha]\cap\{x\in\alpha
X:x\le\alpha\}$ and $\mathcal{C}''=\mathcal{C}[\alpha]\cap\{x\in\alpha
X:x\ge\alpha\}$ are both connected subsets of $X$ with $\alpha\in
\mathcal{C}'\cap\overline{(X\setminus\mathcal{C}')}$ and
$\alpha\in\mathcal{C}''\cap\overline{(X\setminus\mathcal{C}'')}$
(consider that $X\setminus\mathcal{C}'\supset R$ and
$X\setminus\mathcal{C}''\supset L$), so that $\mathcal{C}'=\{\alpha\}$
and $\mathcal{C}''=\{\alpha\}$, whence also
$\mathcal{C}[\alpha]=\{\alpha\}$.
Then, $\alpha X$ has a clopen base
at $\alpha$. Indeed, let $\ell\in L$ and $r\in R$. Since
$\mathcal{C}[\alpha]$ is also the quasi-component of the point
$\alpha$, there are clopen neighbourhoods $U_\ell,U_r$ of $\alpha$
such that $\ell\not\in U_\ell$ and $r\not\in U_r$. Then,
\[
U=\{x\in U_\ell\cap U_r: \ell <x<r\}=\{x\in U_\ell\cap U_r: \ell\le
x\le r\}
\]
is a clopen neighbourhood of $\alpha$ with $U\subset
\{x\in X: \ell< x< r\}$.

That is, $\alpha X$ has a clopen base at $\alpha$ and it is first
countable at this point. Then, let $\{U_n:n<\omega\}$ be a decreasing
clopen base at $\alpha$, with $U_0=\alpha X$. Next, for every point
$x\in X$, let $n(x)=\max\{n:x\in U_n\}$ and, for convenience,
$n(\alpha)=\omega$. Finally, we may define a linear order ``$\prec$''
on $\alpha X$ by putting $x\prec y$ if and only if
\[
\begin{array}{llll}
\hbox{either}&
n(x) < n(y)&
\mbox{or}&
n(x) = n(y)\ \mbox{and}\ x< y.
\end{array}
\]
Since $\{U_n:n<\omega\}$ is a decreasing clopen base at $\alpha$, the 
order ``$\prec$'' is compatible with the topology of $\alpha X$. 
It is clear that, with respect to ``$\prec$'', $\alpha$ is the last 
element of $X$. 
This completes the proof.

\providecommand{\bysame}{\leavevmode\hbox to3em{\hrulefill}\thinspace}
\providecommand{\MR}{\relax\ifhmode\unskip\space\fi MR }
\providecommand{\MRhref}[2]{%
  \href{http://www.ams.org/mathscinet-getitem?mr=#1}{#2}
}
\providecommand{\href}[2]{#2}

\end{document}